\documentclass[seceqn,number,dvips]{arxbj}
\usepackage{cursive}

\aid{0}
\volume{13}
\issue{3}
\pubyear{2007}
\firstpage{653}
\lastpage{671}
\doi{10.3150/07-BEJ6025}

\makeatletter

\newtheorem{Thm}{Theorem}
\newtheorem{Lem}{Lemma}
\newtheorem{Pro}{Proposition}
\newremark{rems}{Remarks}
\newremark{ex}{Example}

\makeatother

\begin{document}
\begin{frontmatter}

\title{High-resolution product quantization for Gaussian processes
under sup-norm distortion}
\runtitle{High-resolution product quantization for Gaussian processes}

\begin{aug}
\author[a]{\fnms{Harald} \snm{Luschgy}\ead[label=e1]{luschgy@uni-trier.de}\thanksref[1]{a1}}
\and
\author[b]{\fnms{Gilles} \snm{Pag\`es}\corref{}\ead[label=e2]{gpa@ccr.jussieu.fr}\thanksref[2]{a2}}
\runauthor{H. Luschgy and G. Pag\`es}
\pdfauthor{H. Luschgy, G. Pages}
\address[a]{Universit\"at Trier, FB IV-Mathematik, D-54286 Trier,
Germany. \printead{e1}}
\address[b]{Laboratoire de Probabilit\'es et Mod\`eles al\'eatoires,
UMR 7599, Universit\'e Paris 6, case 188, 4,
pl. Jussieu, F-75252 Paris Cedex 5, France. \printead{e2}}
\end{aug}

\received{\smonth{1} \syear{2005}}
\revised{\smonth{1} \syear{2007}}

%
\begin{abstract}
We derive high-resolution upper bounds for optimal product quantization
of pathwise continuous Gaussian
processes with respect to the supremum norm on $[0,T]^d$. Moreover, we
describe a product
quantization design which attains this bound. This is achieved under
very general assumptions on
random series expansions of the process. It turns out that product
quantization is asymptotically
only slightly worse than optimal functional quantization. The results
are applied to
fractional Brownian sheets and the Ornstein--Uhlenbeck process.

\end{abstract}

%

\end{frontmatter}

\section{Introduction}\label{s1}

In this paper, we investigate the functional quantization problem for
pathwise continuous Gaussian processes
$X = (X_t)_{t\in I}, I = [0,T]^d$, where the path space
$E = \mathcal{C} (I)$ is endowed with the supremum norm. For any real
separable space $(E,\| \cdot\|)$ and
$r \in(0,\infty)$, optimal quantization means the best approximation
in $L^r_E(\mathbb{P})$ of a random
vector $X\dvtx(\Omega,\mathcal{A}, \mathbb{P})\to E$ by random vectors
$\widehat
X\dvtx(\Omega,\mathcal{A}, \mathbb{P})\to E$ taking
finitely many values in $E$. If $N \in\mathbb{N}$, $\operatorname{card}(\widehat
X(\Omega))\le N$, then $\widehat X$ is called
$N$-quantization. This leads to the minimal level $N$\textit{-quantization}
error defined by
%
\begin{equation}
e_{N,r}(X,E) := \inf\{(\mathbb{E}\| X-\widehat X \|^r )^{1/r} \dvtx
\widehat X N\mbox{-quantization of } X\},
\end{equation}
provided $X \in L^r_E(\mathbb{P})$. When $E=\mathbb{R}^d$, this
problem is known as
\textit{optimal vector quantization} and has been extensively
investigated since the early 1950s, with some applications to signal
processing and transmission (see Gersho and Gray \cite{GEGR}) and to
model-based
clustering in statistics (see e.g., Tarpey \cite{TA}). Beyond these
classical applications, optimal quantization has been used as a
space discretization device to solve nonlinear problems, such as those
arising in optimal stopping theory (American-style option pricing,
reflected BSDE, Bally and Pag\`es \cite{BAPA}), nonlinear filtering
(Pag\`es and Pham \cite{PAPH}),
forward-backward SDE (see Delarue and Menozzi \cite{DEME}) and SPDE
(see Gobet \textit{et al.} \cite{GOPAPHPR}). The mathematical
foundations are treated
in Graf and Luschgy \cite{GRAF1}. Much attention has been paid to the
infinite-dimensional case. This is the so-called \textit{functional
quantization} of stochastic processes: the aim is to quantize some
processes viewed as random vectors taking values in their
path spaces. Recently, a first application of functional
quantization to statistical clustering of functional data has been
investigated (see Tarpey and Kinateder \cite{TAKI} and Tarpey \textit
{et al.} \cite{TAPEOG}). The simplest
application of functional quantization as a numerical method is to
use it as an alternative to Monte Carlo simulation, using the
quadrature formula
\[
\mathbb{E}(F(X)) \approx\mathbb{E}(F(\widehat X))= \sum_{a\in
\alpha} F(a) \mathbb{P}
(\widehat
X=a),\qquad\mbox{where } \alpha= \widehat X(\Omega),
\]
for sufficiently regular functionals $F\dvtx E\to\mathbb{R}$. If
$\widehat X$
is an
$L^r$-optimal $N$-quantization and
$F$ is Lipschitz continuous, then the induced error is bounded by
$[F]_{\mathrm{Lip}}e_{N,r}(X,E)$,
$r\ge1$. Some numerical applications are being developed for the pricing
of path-dependent options (such as regular Asian options) in various
models using $E= L^{2}([0,T],dt)$ (Black and Scholes, Heston,
see Pag\`es and Printems \cite{PAPR2}, Wilbertz \cite{WI}). However,
many important functionals
of processes, like those related to barrier options or to options
on maximum, are only continuous with respect to the sup-norm on
$E = \mathcal{C}([0,T])$.

Let us now describe what we will call the \textit{product
quantization} scheme. Let $X$ be a centered
$E$-valued Gaussian random vector. Let
$\xi_1, \xi_2, \ldots$ be i.i.d. $\mathcal{N}(0,1)$-distributed random
variables and let $(f_j)_{j \geq1}$ be
a sequence in $E$ such that $\sum^\infty_{j=1} \xi_j f_j$
converges a.s. in $E$ and
%
\begin{equation}
X \stackrel{d}{=} \sum^\infty_{j=1} \xi_j f_j .
\end{equation}
Let us call such a sequence \textit{admissible} for $X$. For background
on expansions for Gaussian random vectors, the reader is referred
to Bogachev \cite{BOGA} and Ledoux and Talagrand \cite{LEDO}. One
checks that $(f_j)_{j\ge1}$
is admissible for $X$ if and only if $(f_j)_{j\ge1}$ is a
normalized tight frame in the reproducing kernel Hilbert space
(Cameron--Martin space) $H=H_{_X}$, that is, $\{f_j, j\ge1\}\subset
H$ and $\sum_{j\ge1} (f_j,h)^{2}_{_H} = \|h\|_{_H}^{2}$ for all
$h \in H$ (see Luschgy and Pag\`es \cite{LUPASpec}). Then a sufficient
(but not
necessary) condition is that $(f_j)_{j\ge1}$ is an orthonormal basis
of $H_{_X}$.

For $m, N_1,\ldots, N_m\in\mathbb{N}$ with $\prod^m_{j=1} N_j \leq N$, let
$\widehat\xi_j$ be an $L^r$-optimal $N_j$-quantization for $\xi_j$,
that is,
$(\mathbb{E}|\xi_j-\widehat\xi_j|^r)^{1/r}= e_{N_j,r}(\xi
_j,\mathbb{R}).$ An
$L^r$-product $N$-quantization of $X$ with respect to $(f_j)_{j \geq
1}$ is then defined by
%
\begin{equation}\label{1.3}
\widehat X :=\widehat X^{(N_1,\ldots,N_m)}:= \sum^m_{j=1} \hat
{\xi}_j f_j
\end{equation}
and the quantization error induced by $\widehat X$ is
\[
( \mathbb{E}\|X - \widehat X \|^r )^{1/r} .
\]
Note that if $\alpha_j= \widehat\xi_j(\Omega)$, then the codebook
$\alpha=
\widehat X(\Omega)$ of $\widehat X$
satisfies $\alpha=\{\sum_{j=1}^ma_jf_j \dvtx a \in\prod
_{j=1}^m\alpha_j\}
$ and
\[
\biggl( \mathbb{E}\min_{a \in\alpha} \| X - a \|^r \biggr)^{1/r}
\leq(
\mathbb{E}\|X -
\widehat X \|^r )^{1/r} .
\]
The minimal $N$th product quantization error is then defined by
%
\begin{eqnarray}
e^{(\mathrm{prod})}_{N,r} (X,E) &:=& \inf\{ (\mathbb{E}\| X - \widehat X
\|^r )^{1/r} \dvtx(f_j)_{j\ge1} \in E^{\mathbb{N}} \mbox{
admissible for }
X, \nonumber\\
%
&& \hspace*{42pt}\qquad\widehat X L^r\mbox{-product } N\mbox{-quantization
w.r.t. } (f_j) \}.
\end{eqnarray}
Clearly, we have
%
\begin{equation}
e_{N,r}(X,E) \leq e^{(\mathrm{prod})} _{ N,r}(X,E).
\end{equation}

We address the issue of high-resolution product quantization in $E =
\mathcal{C}(I)$ under the sup-norm, which concerns the performance of
$\widehat X=\widehat X^{(N_1,\ldots,N_m)}$ under a suitable choice
of the marginal quantization levels $N_j$ and the behaviour of
$e^{(\mathrm{prod})}_{ N,r}(X, \mathcal{C}(I))$ as $N \rightarrow\infty$.
For a broad class of Gaussian processes, we derive high-resolution
upper estimates for $e^{(\mathrm{prod})}_{ N,r}(X, \mathcal{C}(I))$.
Furthermore, we describe a product quantization design $\widehat X$
which attains this bound. Combining these estimates with precise
high-resolution formulas for $e_{ N,r}(X, \mathcal{C} (I))$ (see
Dereich \textit{et al.} \cite{DEREI1}, Dereich and Scheutzow \cite{DEREI2},
Graf \textit{et al.} \cite{GRAF2}), one may typically conclude that
\[
e^{(\mathrm{prod})} _{ N,r}(X, \mathcal{C} (I))= O \bigl(( \log\log N)^c e_{
N,r} (X, \mathcal{C}(I))\bigr),
\]
for some suitable constant $c > 0$. This suggests that the asymptotic
quality of product quantization, which is based on easy computations,
is only slightly worse than optimal quantization. 
The optimality of this rate for product quantization rate remains open,
although one may reasonably guess that
it is optimal.

The paper is organized as follows. In Section \ref{s2}, we derive
high-resolution upper estimates for $e^{(\mathrm{prod})}_{ N,r}(X, \mathcal
{C}(I))$ under very general assumptions on expansions. Section \ref{s3}
contains a collection of examples, including fractional Brownian sheets,
Riemann--Liouville processes and the Ornstein--Uhlenbeck process.

It is convenient to use the symbols $\sim$ and $\approx$, where $a_n
\sim b_n$ means $a_n/b_n \rightarrow1$ and
$a_n \approx b_n$ means $a_n = O(b_n)$ and $a_n = \Omega(b_n)$.
Throughout, all logarithms are natural logarithms and
$[x]$ denotes the integer part of the real number
$x$.

\section{High-resolution product quantization}\label{s2}

We investigate high-resolution product functional quantization of
centered continuous Gaussian processes
$X = (X_t)_{t \in I}$ on $I = [0,T]^d$ in the space $E = \mathcal{C}(I)$
equipped with the sup-norm
$\| x \| = \sup_{t \in I} | x (t) | $. Let
%
\[
e^{(\mathrm{prod})}_{ N,r} (X):= e^{(\mathrm{prod})}_{ N,r} (X, \mathcal{C}(I)).
\]
The subsequent setting comprises a broad class of processes.

Let $(f_j)_{j \geq1} \in\mathcal{C}(I)^\mathbb{N}$ satisfy the
following assumptions:
\begin{longlist}
\item[(A1)]
$\| f_j \| \leq C_1 j^{-\vartheta} \log(1+j)^\gamma$ for every $j
\geq
1$ with $\vartheta> 1/2, \gamma\geq0$ and $C_1 < \infty$;

\item[(A2)]
$f_j$ is $a$-H{\"o}lder-continuous and $[f_j]_a \leq C_2 j^b$ for every
$j \geq1$ with $a \in(0,1], b \in\mathbb{R}$ and $C_2 < \infty$, where
\[
[f]_a = \sup_{s \not= t} \frac{ | f(s) - f(t) | }{ | s-t |^a }
\]
(and $| t | $ denotes the $l_2$-norm of $t \in\mathbb{R}^d$).
\end{longlist}

In the sequel, finite constants depending only on the parameters $T,
\vartheta, \gamma, a, b, C_1, C_2, d$ and~$r$ are denoted by $C$ and
may differ from one formula to another one. Other dependencies
are explicitly indicated.

First, observe that by (A1), $\sum^\infty_{j=1} f_j (t)^{2} \leq
\sum^\infty_{j=1} \| f_j \|^{2} < \infty$ for every $t \in I$,
so we can define a centered Gaussian process $Y$ by $Y_t :=
\sum^\infty_{j=1} \xi_j f_j (t)$. Using (A1) and (A2), we
have, for \mbox{$\rho\in(0,1]$},
\begin{eqnarray*}
| f_j (s) - f_j (t) | & = & | f_j (s) - f_j (t) |^\rho| f_j (s) - f_j
(t) |^{1-\rho} \\
 & \leq& ( [f_j]_a | s-t |^a )^\rho(2 \| f_j \| )^{1-\rho} \\
 & \leq& C_\rho j^{\rho(b+\vartheta) - \vartheta} \log
(1+j)^{\gamma(1-\rho)} | s-t |^{a \rho}
\end{eqnarray*}
and hence
\[
\sum^\infty_{j=1} [f_j]^{2}_{a \rho} < \infty\qquad \mbox{for
every } \rho< \frac{\vartheta- 1/2}{(b+\vartheta)_+} .
\]
This yields
%
\begin{equation}
\mathbb{E}| Y_s - Y_t |^{2} = \sum^{\infty}_{j=1} | f_j (s) - f_j
(t) |^{2}
\leq\Biggl( \sum^\infty_{j=1} [f_j]^{2}_{a \rho} \Biggr) | s-t |^{2 a \rho}
\end{equation}
and using the Gaussian feature of $Y$, we obtain from the Kolmogorov
criterion that $Y$ has a continuous
modification $X$. Consequently, $(f_j)$ is admissible for $X$ and
%
\begin{equation}
X = \sum^\infty_{j=1} \xi_j f_j \qquad\mbox{a.s.}
\end{equation}
For $r \in[1, \infty)$, the quantization error induced by the
$L^r$-product $N$-quantization $\widehat X:=\widehat
X^{(N_1,\ldots,N_m)}$ (see (\ref{1.3})) satisfies
%
\begin{eqnarray*}
\nonumber(\mathbb{E}\| X - \widehat X\|^r)^{1/r} & = & \| X -
\widehat X\|
_{L^r_E(\mathbb{P})} \\
\nonumber& \leq& \Biggl\| \sum^m_{j=1} ( \xi_j - \hat{\xi}_j) f_j \Biggr\|
_{L^r_E(\mathbb{P})} + \Biggl\|
\sum_{j \geq m+1} \xi_j f_j \Biggr\|_{L^r_E(\mathbb{P})} \\
\nonumber& \leq& \sum^m_{j=1} \| \xi_j - \hat{\xi}_j \|
_{L^r(\mathbb{P}
)} \| f_j \| + \Biggl\| \sum_{j \geq m+1} \xi_j f_j \Biggr\|_{L^r_E(\mathbb{P})}\\
\end{eqnarray*}
so that
\begin{equation}
\label{rge1}
( \mathbb{E}\| X - \widehat X \|
^r)^{1/r} \leq \sum^m_{j=1} \| f_j \| e_{N_j, r}
(\mathcal{N}(0,1)) + \Biggl( \mathbb{E}\Biggl\| \sum_{j \geq m+1} \xi_j f_j \Biggr\|
^r\Biggr)^{1/r} .
\end{equation}
For $r \in(0,1)$, we have
%
\begin{equation}\label{rle1}
(\mathbb{E}\| X - \widehat X\|^r)^{1/r}\le\mathbb{E}\| X -
\widehat X\| \leq
\sum^m_{j=1} \| f_j \| \mathbb{E}|\xi_j-\widehat\xi_j|+ \mathbb
{E}\Biggl\|
\sum_{j \geq m+1} \xi_j f_j\Biggr\|.
\end{equation}

Let us now consider the truncation error.
\begin{Thm}\label{thm1}
Assume that $(f_j)_{j \geq1} \in\mathcal{C} (I)^\mathbb{N}$ satisfies
(\textup{A1})--(\textup{A2}). Then, for every $n \geq2$ and \mbox{$r \in(0, \infty)$},
\[
\Biggl(\mathbb{E}\Biggl\| \sum_{j \geq n} \xi_j f_j \Biggr\|^r \Biggr)^{1/r} \leq
\frac{C(\log n)^{\gamma+ 1/2}}{n^{\vartheta- 1/2}}
\]
and
\[
\Biggl( \mathbb{E}\Biggl\| \sum_{j \geq n} \xi_j f_j \Biggr\|^r \Biggr)^{1/r} \leq
\frac{C(\log n)^\gamma}{n^{\vartheta- 1/2}},\qquad \mbox{if } b +
\vartheta\leq0.
\]
\end{Thm}
\begin{pf}
By equivalence of Gaussian moments,
%
\begin{equation}
\Biggl( \mathbb{E}\Biggl\| \sum_{j \geq n} \xi_j f_j \Biggr\| ^r\Biggr)^{1/r} \leq D \mathbb
{E}\Biggl\|
\sum_{j \geq n} \xi_j f_j \Biggr\|
\end{equation}
for some constant $D$ depending on $r$ (cf. Ledoux and Talagrand \cite{LEDO}, Corollary
3.2). The upper estimate for $\mathbb{E}\| \sum_{j \geq n} \xi_j f_j
\| $ is based on corresponding estimates for finite blocks of
exponentially increasing length. For $m \geq1$, set
\[
Z = Z^{(m)} := \sum^{2^m}_{j=2^{m-1} +1} \xi_j f_j .
\]
For a given $N \geq1$, consider the grid $G_{ N} = \{ \frac{
(2i-1)T}{2N} \dvtx i = 1, \ldots, N \}^d$. Then
\[
\| Z \| \leq\sup_{t \in G_N} | Z_t | + \sup_{| s-t| \leq C N^{-1}}
| Z_s - Z_t | .
\]
It follows from the Gaussian maximal inequality that
\[
\mathbb{E}\sup_{t \in G_N} | Z_t | \leq C \sqrt{\log(1+N^d)} \sup
_{t \in G_N}
\sqrt{\mathbb{E}Z^{2}_t} .
\]
Using (A1), we have, for every $t \in I$,
\[
\mathbb{E}Z^{2}_t \leq\sum^{2^m}_{j=2^{m-1}+1} \| f_j \|^{2} \leq C
\sum^{2^m}_{j=2^{m-1}+1}
j^{-2 \vartheta}\log(1+j)^{2 \gamma} \leq C 2^{m(1-2 \vartheta)} m^{2
\gamma}
\]
so that
\[
\mathbb{E}\sup_{t \in G_{N}} | Z_t | \leq C \sqrt{\log(1+N)} 2^{-m
(\vartheta
-1/2)} m^\gamma.
\]
Moreover, using (A2), we have, for $| s-t | \leq C N^{-1}$,
\begin{eqnarray*}
| Z_s - Z_t | & \leq& \sum^{2^m}_{j=2^{m-1}+1} | \xi_j |
| f_j (s) - f_j (t) | \\
& \leq& C | s-t |^a \sum^{2^m}_{j=2^{m-1}+1} | \xi_j | [f_j]_a
\\
& \leq& C N^{-a} \sum^{2^m}_{j=2^{m-1}+1} | \xi_j | j^b
\end{eqnarray*}
and hence
\[
\mathbb{E}\sup_{| s-t | \leq C
N^{-1}} | Z_s - Z_t | \leq C N^{-a}
\sum^{2^m}_{j=2^{m-1}+1} j^b
\leq C N^{-a} 2^{m(1+b)} .
\]
Thus we have established the estimate
%
\begin{equation}\label{(2.6)}
\mathbb{E}\bigl\| Z^{(m)} \bigr\| \leq C \bigl( \sqrt{\log(1+N)} 2^{-m(\vartheta-1/2)}
m^\gamma+ N^{-a} 2^{m (1+b)}\bigr).
\end{equation}
As concerns the choice of $N$, set $N := [2^{u m}] +1$, with $u \in
(0,\infty)$ satisfying
$1+b-au \leq\frac{1}{2} - \vartheta$.
Equation (\ref{(2.6)}) then becomes
%
\begin{equation}
\mathbb{E}\bigl\| Z^{(m)} \bigr\| \leq C 2^{-m(\vartheta-1/2)} m^{\gamma+1/2} .
\end{equation}

We note that in the case $b + \vartheta\leq- 1/2$, we may choose $N =
1$ and thereby obtain a power reduction from
$m^{\gamma+ 1/2}$ to $m^\gamma$. This can be improved. In fact, we have
\begin{eqnarray*}
\mathbb{E}| Z_s - Z_t |^{2} & = & \sum^{2^m}_{j=2^{m-1}+1} | f_j (s) -
f_j(t) |^{2} \\
& \leq& C | s-t |^{2a} \sum^{2^m}_{j-2^{m-1}+1} j^{2b} \leq C |
s-t |^{2a}
2^{m(1+2b)}
\end{eqnarray*}
so that
\[
d_Z (s,t) := ( \mathbb{E}| Z_s - Z_t |^{2} )^{1/2} \leq C | s-t |^a
2^{m(b+1/2)} .
\]
If $N(\varepsilon, d_Z)$ denotes the covering numbers of $I$ with
respect to the intrinsic semi-metric $d_Z$,
then, by chaining,
\[
\mathbb{E}\sup_{| s-t| \leq C N^{-1}} | Z_s - Z_t | \leq\mathbb
{E}\sup_{d_Z (s,t)
\leq\delta} |
Z_s - Z_t | \leq C \int^\delta_0 \sqrt{\log N(\varepsilon, d_Z)} \,\mathrm{d}
\varepsilon,
\]
where $ \delta:= C N^{-a} 2^{m (b+1/2)}$ (cf. Van der Waart and
Wellner \cite{WAA}, page 101).
Since
\[
N(\varepsilon, d_Z) \leq C \biggl( \frac{2^{m(b+1/2)}}{\varepsilon}\biggr)^{d/a}
,\qquad 0<\varepsilon\leq\varepsilon_0,
\]
and $\int^1_0 \sqrt{\log(1/x)} \,\mathrm{d}x < +\infty$, we obtain, for sufficiently large $N$,
\[
\int^\delta_0 \sqrt{\log N(\varepsilon, d_Z)} \,\mathrm{d} \varepsilon\leq C
2^{m(b+1/2)} \int^1_0
\sqrt{\log(1/x)} \,\mathrm{d}x \leq C 2^{m(b+1/2)} .
\]
Consequently,
%
\begin{eqnarray}
\mathbb{E}\bigl\| Z^{(m)} \bigr\| & \leq& C \bigl( \sqrt{\log(1+N)}
2^{-m(\vartheta-1/2)}
m^\gamma+ 2^{m(b+1/2)}\bigr) \nonumber\\
& \leq& C 2^{-m(\vartheta-1/2)} m^\gamma\qquad\mbox{if }
b+\vartheta\leq0 .
\end{eqnarray}

We now complete the proof. For $n \geq2$, choose $m = m(n) \geq1$
such that $2^{m-1} < n \leq2^m$. Then
\[
\Biggl\| \sum_{j \geq n} \xi_j f_j \Biggr\| \leq\sum_{j \geq m+1} \bigl\|
Z^{(j)} \bigr\|
+ \Biggl\| \sum^{2^m}_{j=n} \xi_j f_j \Biggr\| .
\]
Since $\mathbb{E}\| \sum_{n\le j\le2^m} \xi_j f_j \| \leq\mathbb
{E}\| Z^{(m)}
\| $ by the Anderson inequality (cf. Bogachev \cite{BOGA}, Corollary
3.3.7), we
deduce from equation (2.7) that
\[
\mathbb{E}\Biggl\| \sum_{j \geq n} \xi_j f_j \Biggr\| \leq C \sum_{j \geq m}
\frac{j^{\gamma+ 1/2}}{2^{j(\vartheta- 1/2)}} \leq
\frac{C m^{\gamma+ 1/2}}{2^{m(\vartheta-1/2)}}
\leq\frac{C ( \log n)^{\gamma+ 1/2}}{n^{\vartheta- 1/2}} .
\]
If $b + \vartheta\leq0$, then it follows from (2.8) that
\[
\mathbb{E}\Biggl\| \sum_{j \geq n} \xi_j f_j \Biggr\| \leq\frac{C ( \log
n)^\gamma
}{n^{\vartheta- 1/2}} .
\]
Combining these estimates with (2.5) yields the assertion.
\vadjust{\goodbreak}\end{pf}

\begin{rems*}
\begin{longlist}
\item[$\bullet$] The rate for the truncation error
depends only on $\vartheta$ and $\gamma$, that is, on the decay of the size
of functions $f_j$ (provided $b + \vartheta> 0$). The occurrence of
expansions with $b + \vartheta\leq0$ seems to be a rare event and
otherwise $b$ plays no role (see the subsequent example). The case
$\gamma= 0$ typically corresponds to one-parameter processes with $I =
[0,T]$.

\item[$\bullet$] The $e^{(\mathrm{prod})}_{N,r}$-problem comprises the
optimization of admissible sequences and, in view of (\ref{rge1})
and (\ref{rle1}), is thus related to the $l$-numbers of $X$ defined by
%
\begin{equation}
l_{n,r} (X) = l_{n,r} (X, \mathcal{C} (I)) := \inf\Biggl\{ \Biggl( \mathbb{E} \Biggl\|
\sum_{j \geq n} \xi_j g_j  \Biggr\|^r \Biggr)^{1/r} \dvtx (g_j)
\mbox{ admissible for } X \mbox{ in } \mathcal{C}(I) \Biggr\} .
\end{equation}
\end{longlist}
\end{rems*}
Rate-optimal solutions of the $l_{n,r}$-problem, in the sense of
$l_{n,r} (X) \approx( \mathbb{E}\| \sum_{j \geq n} \xi_j g_j
\|^r)^{1/r}$ as $n \rightarrow\infty$, have recently been investigated
(see K{\"u}hn and Linde \cite{KUHN}, Dzhaparidze and van Zanten
\cite{DZHA,DZHA1,DZHA2}, Ayache and Taqqu \cite{AYA}). Admissible
sequences of type (A1) and (A2) seem to be promising candidates.

\begin{ex}[{\normalfont(\textit{Weierstrass processes})}]
Let
\[
f_j (t) =j^{- \vartheta} \sin(j^{b+ \vartheta} t),\qquad j \geq1,
\vartheta> 1/2, b \in\mathbb{R}, t \in[0,T].
\]
Then $\| f_j \| \leq j^{-\vartheta}$ and $[f_j]_1 = j^b$. Since $f_j
(0) = 0$, we also have
$\| f_j \| \leq T j^b$, so (A1) and (A2) are satisfied, with
$\tilde{\vartheta} = \max\{ \vartheta, -b\}$ and $a = 1$. The
covariance function of $X = \sum^\infty_{j=1} \xi_j f_j$ is
given by
\[
\mathbb{E}X_s X_t = \sum_{j\ge1} j^{-2 \vartheta} \sin
(j^{b+\vartheta} s)
\sin(j^{b+\vartheta} t).
\]
Now, in the ``Weierstrass case'' $b + \vartheta> 0$, we obtain, from
Theorem \ref{thm1},
\[
\Biggl( \mathbb{E}\Biggl\| \sum_{j \geq n} \xi_j f_j \Biggr\| ^r\Biggr)^{1/r} \leq
\frac{C \sqrt{\log n}}{n^{\vartheta- 1/2}},
\]
while in the ``non-Weierstrass case,'' $b + \vartheta\leq0$ appears
the better rate:
\[
\Biggl( \mathbb{E}\Biggl\| \sum_{j \geq n} \xi_j f_j \Biggr\|^r\Biggr)^{1/r} \leq
\frac{C}{n^{-b-1/2}}.
\]

We pass to the minimal product quantization error $e^{(\mathrm{prod})}_{ N,r}(X)$.
\end{ex}

\begin{Thm}\label{thm2}
Assume that $X$ admits an admissible set $(f_j)_{j \geq1 }$ in
$\mathcal{C}(I)$ satisfying (\textup{A1}) and (\textup{A2}). We then have, for every
$N \geq3$
and $r \in(0, \infty)$,
%
\begin{equation}
e^{(\mathrm{prod})}_{ N,r} (X) \leq\frac{C (\log\log N)^{\vartheta+
\gamma
} }{(\log N)^{\vartheta- 1/2}}
\end{equation}
and
\[
e^{(\mathrm{prod})}_{ N,r}(X)\leq\frac{C(\log\log N)^{\vartheta+ \gamma-
1/2} }{ (\log N)^{\vartheta-1/2} } \qquad
\mbox{if } b + \vartheta\leq0 .
\]
Furthermore, the $L^r$-product $N$-quantization $\widehat X$ with
respect to $(f_j)$, with tuning parameters defined in \textup{(2.11)} and \textup{(2.15)}
below, achieves these rates.
\end{Thm}

\begin{pf}
Let $r \in[1, \infty)$ and set
$\nu_j := j^{-\vartheta}_0 \log(1+j_0)^\gamma$
if $j < j_0 := [e^{\gamma/\vartheta}]$ and $\nu_j := j^{-\vartheta}
\log(1+j)^\gamma$ if
$j \geq j_0$. The sequence $(\nu_j)_j$ is then decreasing. Since
\[
\lim_{k \to\infty} k e_{k,r} (\mathcal{N}(0,1),\mathbb{R}) \mbox
{ exists in }
(0, \infty)
\]
(cf. Graf and Luschgy \cite{GRAF1}), we deduce from (\ref{rge1}),
(A1) and Theorem \ref{thm1} the estimate
\[
( \mathbb{E}\| X - \widehat X\|^r)^{1/r} \leq C \Biggl( \sum^{m}_{j=1} \nu_j
N^{-1}_j +
\frac{\log(1+m)^{\gamma+1/2}}{m^{\vartheta- 1/2}}\Biggr),
\]
for every $m, N_1 , \ldots, N_m \in\mathbb{N}$ with $\prod^m_{j=1} N_j
\leq N$.
(The case $b+\vartheta\leq0$ is treated analogously.) Consequently,
%
\begin{eqnarray}
e^{(\mathrm{prod})}_{ N,r} (X) &\leq& C \inf\Biggl\{ \sum^{m}_{j=1} \nu_j
N^{-1}_j +
\frac{\log(1+m)^{\gamma+1/2}}{m^{\vartheta- 1/2}} \dvtx
m, N_1 , \ldots, N_m \in\mathbb{N}, \nonumber\\
&&\hspace*{189pt}
\prod^m_{j=1} N_j \leq N \Biggr\} .
\end{eqnarray}
For a given $N \in\mathbb{N}$, we may first optimize the integer bit
allocation given by the $N_j$'s for fixed $m$
and then optimize $m$. To this end, note that the continuous allocation
problem reads
\[
\inf\Biggl\{ \sum^{m}_{j=1} \nu_j y^{-1}_j \dvtx y_j > 0, \prod^m_{j=1} y_j
\leq N \Biggr\} = \sum^{m}_{j=1}
\nu_j z^{-1}_j = N^{-1/m} m \Biggl( \prod^m_{j=1} \nu_j\Biggr)^{1/m},
\]
where
\[
z_j = N^{1/m} \nu_j \Biggl( \prod^m_{k=1} \nu_k\Biggr)^{-1/m}
\]
and $z_1 \geq\cdots\geq z_m$. One can produce an (approximate)
integer solution by setting
%
\begin{equation}
N_j = [z_j] = \Biggl[N^{1/m} \nu_j \Biggl( \prod^m_{k=1} \nu_k \Biggr)^{-1/m} \Biggr] ,\qquad j \in
\{1, \ldots, m \},
\end{equation}
provided $z_m \geq1$. Then
\[
\sum^{m}_{j=1} \nu_j N^{-1}_j \leq2 m N^{-1/m} \Biggl( \prod^m_{j=1} \nu
_j\Biggr)^{1/m} \leq C m N^{-1/m} \nu_m.
\]
Since the constraint on $m$ reads $m \in I(N)$ with
%
\begin{equation}
I (N) := \Biggl\{ m \in\mathbb{N}\dvtx N^{1/m} \nu_m \Biggl( \prod^m_{j=1} \nu_j \Biggr)^{-1/m}
\geq 1 \Biggr\},
\end{equation}
we arrive at
%
\begin{equation}
e^{(\mathrm{prod})} _{ N,r} (X)\leq C \inf_{m \in I(N)} \biggl( \frac{N^{-1/m}
\log(1+m)^\gamma}{ m^{\vartheta-1} } +
\frac{\log(1+m)^{\gamma+1/2}}{m^{\vartheta-1/2}} \biggr),
\end{equation}
for every $N \in\mathbb{N}$. We check that $I(N)$ is finite, $I(N)
= \{ 1,
\ldots, m^{*} (N) \}$, $m^{*}(N)$
increases to infinity and
%
\begin{equation}
m^{*} (N) \sim\frac{\log N}{\vartheta} \qquad \mbox{as } N \rightarrow
\infty.
\end{equation}
Finally, let
%
\begin{equation}
m  =  m(N) \in I(N),\qquad \mbox{with } m(N) \leq\frac{2 \log N}{\log\log
N} \qquad \mbox{for }
N \geq3
\end{equation}
such that
\[
m(N) \sim \frac{2 \log N}{\log\log N}
\qquad\mbox{as } N \rightarrow\infty.
\]
This is possible in view of (2.14). Using (\ref{rle1}), the case $r
\in(0,1)$ follows from $r=1$ since the $L^r$-optimal
$N_j$-quantizations $\widehat\xi_j$ satisfy $\mathbb{E}|\xi
_j-\widehat
\xi_j| \le C N_j^{-1}$, $j\ge1$; see  Graf \textit{et al.} \cite{GRAF3}.
\end{pf}

We may reasonably conjecture that for many specific processes, the
above rate is the true one. This would imply that product
quantization achieves the optimal rate for quantization, namely the
rate of convergence to zero of $e_{ N,r}(X) := e_{ N,r}(X,
\mathcal{C}(I))$, only up to a $\log\log N$ term in formula
(2.16). This is in contrast to the Hilbert space setting, where the
optimal rate is attained by product quantization (cf. Luschgy and Pag\`es \cite{LUS2}).
To be precise, we summarize the results on $e_{ N,r }(X)$ in the
present setting.
\begin{Pro}\label{pro1}
\textup{(a)} Assume that $X$ admits an admissible sequence in $\mathcal{C}(I)$
satisfying (\textup{A1}) and~(\textup{A2}). Then
%
\begin{equation}
e_{ N,r} (X) = O \biggl( \frac{ ( \log\log N)^{\gamma+ 1/2}}{( \log
N)^{\vartheta-1/2}} \biggr)
\end{equation}
and
%
\begin{equation}
e_{ N,r}(X) = O \biggl( \frac{ ( \log\log N)^{\gamma}}{( \log
N)^{\vartheta-1/2}} \biggr), \qquad
\mbox{if } b + \vartheta\leq0 .
\end{equation}

\textup{(b)} Assume that $X$ admits an admissible sequence satisfying (\textup{A1}).
Let $\mu$ be a finite Borel measure on $I$ and let
$V \dvtx \mathcal{C}(I) \rightarrow L^{2} (I, \mu)$ denote the natural
embedding. Then
\[
e_{ N,r} \bigl(V(X), L^{2} (\mu)\bigr) = O\biggl( \frac{ ( \log\log N)^\gamma
}{(\log N)^{\vartheta-1/2}} \biggr)
\]
and
\[
e^{(\mathrm{prod})}_{N,2} \bigl(V(X), L^{2}(\mu)\bigr) = O \biggl( \frac{ ( \log\log
N)^\gamma}{(\log N)^{\vartheta-1/2}} \biggr).
\]
\end{Pro}
\begin{pf}
(a) The proof is not constructive. We use
Proposition 4.1 in Li and Linde \cite{LI}, which relates $l$-numbers
(see (2.9))
and small ball probabilities (but this relation is not always
sharp). By combining this relation and Theorem \ref{thm1}, we obtain
\begin{eqnarray*}
-\log\bigl(\mathbb{P}(\|X\|\le\varepsilon)\bigr)&=& O\biggl(\varepsilon^{-
1/(\vartheta-1/2)}\biggl(\log\biggl(\frac{1}{\varepsilon}\biggr)\biggr)^{(\gamma+1/2)/(\vartheta-1/2)}\biggr),\\
-\log\bigl(\mathbb{P}(\|X\|\le\varepsilon)\bigr)&=& O\biggl(\varepsilon^{-1/(\vartheta
-1/2)}\biggl(\log\biggl(\frac{1}{\varepsilon}\biggr)\biggr)^{\gamma/(\vartheta
-1/2)}\biggr),\qquad\mbox{if } b+\vartheta\le0
\end{eqnarray*}
as $\varepsilon\to0$. We may then apply a known, precise relationship between
these probabilities and $e_{ N,r}(X)$ (cf. Dereich \textit{et al.}
\cite{DEREI1},
Graf \textit{et al.} \cite{GRAF2}) and this leads to the desired estimate.

(b) Let $(f_j)_{j \geq1}$ be an admissible sequence in $\mathcal
{C}(I)$ for $X$ satisfying (A1) and consider an
$L^{2}$-product $N$-quantization of $V(X)$ based on $(V f_j)_{j\ge1}$,
\[
\widehat{V(X)}^N = \sum^m_{j=1} \hat{\xi}_j V (f_j),
\]
where $\widehat\xi_j$ are $L^{2}$-optimal Voronoi
$N_j$-quantizers; see Luschgy and Pag\`es \cite{LUS1}. Then, using the
independence of
$\xi_j - \hat{\xi}_j$, $j \geq1$, and the stationarity property
$\hat{\xi_j} = \mathbb{E}( \xi_j | \hat{\xi}_j)$ of the quantization
$\hat{\xi}_j$, we have
\begin{eqnarray*}
&& \mathbb{E}\Biggl\| \sum^\infty_{j=1} \xi_j V (f_j)  -  \widehat
{V(X)}^N \Biggr\|
^{2}_{L^{2}(\mu)} \\
&&\quad =  \sum^m_{j=1} \mathbb{E}| \xi_j - \hat{\xi}_j |^{2} \| V f_j
\|^{2}_{L^{2}(\mu)} +
\sum_{j \geq m + 1} \| V f_j \|^{2}_{L^{2}(\mu)} \\
&&\quad \leq C \Biggl( \sum^m_{j=1} N^{-2}_j j^{-2 \vartheta} \log
(1+j)^{2 \gamma} + \sum_{j \geq m + 1}
j^{-2\vartheta} \log(1+j)^{2 \gamma}\Biggr).
\end{eqnarray*}
We then argue along the lines of Luschgy and Pag\`es \cite{LUS1} to conclude that
\[
e^{(\mathrm{prod})}_{N,2} \bigl(V(X), L^{2} (\mu)\bigr) = O \biggl( \frac{ ( \log\log
N)^\gamma
}{ (\log N)^{\vartheta-1/2} } \biggr) .
\]
\upqed
\end{pf}

Sometimes, (2.17) provides the true rate for $e_{ N,r}(X)$ (as for
the two-parameter Brownian sheet),
sometimes it yields the best known upper bound (as for the
$d$-parameter Brownian sheet with $d \geq3$)
and sometimes (2.18) provides the true rate (as for Brownian motion).
The latter typically occurs when the rate of $e_{ N,r }(X)$ and
the ``Hilbert rate'' of
$e_{ N,r}(V(X), L^{2}(\mathrm{d}t))$ coincide (see Section \ref{s3}).
It remains an open question to find conditions for this to happen.

%
%
%

%

\section{Examples}\label{s3}

%
\subsection{Fractional Brownian motions and fractional Brownian
sheets}\label{s3.1}
We consider the Dzaparidze--van Zanten expansion of the \textit{fractional
Brownian motion}\break
$X = (X_t)_{t \in[0,T]}$ with Hurst index $\rho\in(0,1)$ and
covariance function
\[
\mathbb{E}X_s X_t = \tfrac{1}{2} (s^{2\rho} + t^{2\rho} - | s-t
|^{2\rho} ) .
\]
These authors discovered, in Dzhaparidze and van Zanten \cite{DZHA1}, that the sequence
%
\begin{eqnarray}\label{(3.1)}
f^1_j (t) & = & \frac{T^\rho c_\rho\sqrt{2}}{ | J_{1-\rho} (x_j) |
x^{\rho+1}_j } \sin\biggl( \frac{x_jt}{T} \biggr),\qquad\hspace*{27pt} j \geq1, \nonumber \\
f^{2}_j (t) & = & \frac{T^\rho c_\rho\sqrt{2}}{ | J_{-\rho} (y_j) |
y^{\rho+1}_j } \biggl(1 - \cos\biggl( \frac{y_jt}{T}\biggr)\biggr),\qquad j \geq1,
\end{eqnarray}
in $\mathcal{C}([0,T])$ is admissible for $X$, where $J_\nu$ denotes the
Bessel function of the first kind of order~$\nu$, $0 < x_1 < x_2 <
\cdots$ are the positive
zeros of $J_{-\rho} $, $0 < y_1 < y_2 < \cdots$ the positive zeros of
$J_{1 - \rho}$ and
$c^{2}_\rho= \Gamma(1 + 2 \rho) \sin(\curpi\rho)/\curpi$.

Using the asymptotic properties
\[
x_j \sim y_j \sim\curpi j ,\qquad J_{1-\rho} (x_j) \sim J_{-\rho} (y_j) \sim
\frac{\sqrt{2}}{\curpi} j^{-1/2} \qquad\mbox{as } j \rightarrow\infty
\]
(cf. Dzhaparidze and van Zanten \cite{DZHA1}), one observes that a
suitable arrangement of the
functions~(\ref{(3.1)}) (like $f_{2j} = f^1_j$, $f_{2j-1} = f^{2}_j)$
satisfies (A1) and
(A2) with parameters $\vartheta= \rho+ 1/2$, $\gamma= 0$, $a = 1$ and
$b = 1/2 - \rho$. Consequently,
%
\begin{equation}
e^{(\mathrm{prod})}_{ N,r} (FBM) = O \biggl( \frac{ ( \log\log N)^{\rho+ 1/2}
}{ (\log N)^\rho} \biggr),
\end{equation}
while (see Dereich and Scheutzow \cite{DEREI2}, Graf \textit{et al.}
\cite{GRAF2})
%
\begin{equation}
e_{ N,r} (FBM) \approx(\log N)^{-\rho} .
\end{equation}

The tensor products of functions (\ref{(3.1)}) are admissible for the
\textit{fractional Brownian sheet} $X$ over
$[0,T]^d$ with covariance function
\[
\mathbb{E}X_s X_t = \bigl( \tfrac{1}{2}\bigr)^d \prod^d_{i=1} (s^{2 \rho_i}_i + t^{2
\rho
_i}_i - | s_i - t_i |^{2 \rho_i} ) ,
\]
$\rho_i \in(0,1)$, and satisfy conditions (A1) and (A2) with
$\vartheta= \rho+ 1/2$, $\rho= \min_{1 \leq i \leq d} \rho_i$, \mbox{$\gamma= \vartheta(m-1)$}, where
$m = \operatorname{card} \{ i \in\{1 , \ldots, d \} \dvtx \rho_i = \rho\}$, $a =
1$ and $b = \max_{1 \leq i \leq d} (1/2 - \rho_i))_{+}$.
This is a consequence of the following lemma which ensures stability of
conditions (A1) and
(A2) under tensor products.

\begin{Lem}\label{lem1}
For $i \in\{1, \ldots, d \}$, let $(f^i_j)_{j \geq1} \in\mathcal{C}
([0,T])^\mathbb{N}$
satisfy (\textup{A1}) and (\textup{A2}) with parameters $\vartheta_i, \gamma_i, a_i,
b_i$ such that $\gamma_i = 0$. Then a decreasing arrangement of
$(\bigotimes^{d}_{i=1} f^i_{j_i})_{\underline{j} \in\mathbb{N}^d}$
satisfies (\textup{A1}) and (\textup{A2}) with parameters
$\vartheta= \min_{1 \leq i \leq d} \vartheta_i$, $\gamma= \vartheta
(m-1)$, where
$m = \operatorname{card} \{ i \in\{ 1, \ldots, d \}\dvtx \vartheta_i =
\vartheta\}$,
$a = \min_{1 \leq i \leq d} a_i$ and $b = (\max_{1 \leq i \leq d} b_i)_{+}$.
\end{Lem}

\begin{pf}
For $\underline{j} = (j_1 , \ldots, j_d) \in\mathbb
{N}^d$, set
$f_{\underline{j}} = \bigotimes^d_{i=1} f^i_{j_i}$ so that
$f_{\underline{j}} (t) = \prod^d_{j=1} f^i_{j_i} (t_i)$, $t \in[0,T]^d$.
We have
\[
\| f_{\underline{j}} \| \leq\prod^d_{i=1} \| f^i_{j_i} \| \leq C
\prod^d_{i=1} j^{-\vartheta_i}_{i}
\quad\mbox{and}\quad
\bigl| f_{\underline{j}} (s) - f_{\underline{j}} (t) \bigr| \leq C \max_{1
\leq i
\leq d} j^b_i | s-t |^a .
\]
Let $u_{\underline{j}} := \prod^d_{i=1} j^{-\vartheta_i}_{i}$. Choose a
bijective map
$\psi\dvtx \mathbb{N}\rightarrow\mathbb{N}^d$ such that $u_k := u_{\psi
(k)}$ is
decreasing in $k \geq1$. Set
$f_k := f_{\psi(k)}$.
Then
\[
u_k \approx C k^{-\vartheta} (\log k)^{\vartheta(m-1)} \qquad\mbox{as }
k \rightarrow\infty
\]
(cf. Papageorgiou and Wasilkowski \cite{PAPA}, Theorem 2.1). Consequently,
\[
\| f_k \| \leq C k^{- \vartheta} ( \log k)^{\vartheta(m-1)}
\]
and, for $\underline{j} = \psi(k)$,
\[
j_i \leq\prod^d_{i=1} j_i \leq\prod^d_{i=1} j^{\vartheta_i / \vartheta}_i
\leq C k ( \log k)^{-(m-1)} \leq C k ,
\]
hence
\[
| f_k (s) - f_k (t) | \leq C k^b | s-t |^a .
\]
\upqed
\end{pf}

Therefore, by Theorem \ref{thm2} and Proposition \ref{pro1},
%
\begin{equation}\label{eq3.4}
e^{(\mathrm{prod})}_{ N,r }(FBS) = O \biggl( \frac{( \log\log N)^{m(\rho+ 1/2)}
}{ (\log N)^\rho} \biggr)
\end{equation}
and
%
\begin{equation}\label{eq3.5}
e_{ N,r }(FBS) = O \biggl( \frac{( \log\log N)^{m(\rho+ 1/2)-\rho} }{
(\log N)^\rho} \biggr) .
\end{equation}
The Hilbert space setting $E = L^{2}([0,T]^d, \mathrm{d}t)$ provides the lower estimate
%
\begin{equation}\label{eq3.6}
e_{ N,r} (FBS) = \Omega\biggl(
\frac{( \log\log N)^{(m-1)(\rho+ 1/2)} }{ (\log N)^\rho} \biggr)
\end{equation}
(see Luschgy and Pag\`es \cite{LUS1,LUS2}). The true rate of
$e_{ N,r}(FBS)$ is known only for the case $m = 1$, where the
true rate is the ``Hilbert rate'' (\ref{(2.6)}) (see Dereich
\textit{et al.} \cite{DEREI1}),
and for the case $m=2$, where (3.5) is the true rate
(see  Belinsky and Linde \cite{BELI}, Graf \textit{et al.} \cite
{GRAF2}). A reasonable conjecture is
that (3.5) is also the true rate for $m \geq3$.
%
\subsection{Riemann--Liouville and other moving average processes}
For $\psi\in L^{2} ([0,T], \mathrm{d}t)$ and a standard Brownian motion $W$, let
\[
X_t = \int^t_0 \psi(t-s) \,\mathrm{d}W_s,\qquad t \in[0,T],
\]
and assume that $X$ has a pathwise continuous modification.
Since
%
\begin{eqnarray}
\mathbb{E}X_s X_t &=& \int^{s \wedge t}_{0} \psi(s-u) \psi(t-u) \,\mathrm{d}u ,
\nonumber \\
f_j (t) & = & \sqrt{\frac{2}{T}} \int^t_0 \psi(t-s) \cos\biggl( \frac
{\curpi
(j-1/2)s}{T} \biggr) \,\mathrm{d}s \nonumber \\
& = & \sqrt{ \frac{2}{T}} \int^t_0 \psi(s) \cos\biggl( \frac{\curpi
(j-1/2)(t-s) }{T} \biggr) \,\mathrm{d}s,\qquad j \geq1,
\end{eqnarray}
is an admissible sequence for $X$. Observe that (3.7) provides well-defined continuous functions, even for
$\psi\in L^1 ([0,T], \mathrm{d}t)$.
\begin{Lem}\label{lem2}
Let $\psi\in L^1([0,1], \mathrm{d}t)$.

\begin{longlist}
\item[(a)] If $\varphi(t) = \int^t_0 | \psi(s) | \,\mathrm{d}s$ is $\beta$-H{\"o}lder continuous
with $\beta\in(0, 1]$,
then the sequence $(f_j)$ from \textup{(3.7)} satisfies (\textup{A2}) with $a = \beta$
and $b=1$.
In particular, if $\psi\in L^{2} ([0,T], \mathrm{d}t)$, then (\textup{A2}) is
satisfied with $a = 1/2$ and $b = 1$.

\item[(b)] If $\psi$ has finite variation over $[0,T]$, then (\textup{A1}) is
satisfied with $\vartheta= 1$ and
$\gamma= 0$.
\end{longlist}
\end{Lem}
\begin{pf}
Let $\lambda_j = (\curpi(j-1/2)/T)^{-2}$. $(a)$ For $s < t$, we have
\begin{eqnarray*}
f_j(s) - f_j (t) &=& \sqrt{ \frac{2}{T}} \biggl\{ \int^s_0 \psi(u) \bigl( \cos\bigl(
(s-u)/\sqrt{\lambda_j}\bigr) -
\cos\bigl((t-u)/\sqrt{\lambda j} \bigr)\bigr) \,\mathrm{d}u\\
&&\hspace*{105pt}
{}- \int^t_s \psi(u) \cos\bigl( (t-u)/\sqrt{\lambda_j}\bigr)\, \mathrm{d}u \biggr\}
\end{eqnarray*}
so that
\[
| f_j (s) - f_j (t) | \leq\sqrt{\frac{2}{T}} \biggl(
\frac{ | s-t | }{ \sqrt{\lambda_j}}
\| \psi\|_{L^1(\mathrm{d}t)}
+ \int^t_s | \psi(u) | \,\mathrm{d}u \biggr) .
\]
(b) We have
\begin{eqnarray*}
f_j (t) & = & - \sqrt{ 2 \lambda_j /T} \int^t_0 \psi(s) \,\mathrm{d} \bigl(\sin
\bigl((t-s)/\sqrt{\lambda_j} \bigr)\bigr) \\
& = & \sqrt{2 \lambda_j / T} \biggl( \psi(0) \sin\bigl(t/\sqrt{\lambda_j} \bigr) +
\int^t_0 \sin\bigl((t-s)/\sqrt{\lambda_j}\bigr) \,\mathrm{d} \psi(s) \biggr)
\end{eqnarray*}
so that
\[
\| f_j \| \leq\sqrt{2 \lambda_j /T} \bigl( | \psi(0) | + \mbox{Var} (
\psi, [0,T])\bigr) .
\]
\upqed
\end{pf}

This lemma yields a universal upper bound,
\[
e^{(\mathrm{prod})}_{ N,r} (X) = O \biggl( \frac{ \log\log N}{(\log
N)^{1/2}}\biggr),
\]
for functions $\psi$ having finite variation.

In the sequel, we do not concern ourselves with improvements of the parameter $b$ in
(A2) since the condition
$b + \vartheta\leq0$ cannot be achieved in this setting.
\begin{Lem}\label{lem3}
Let $\psi\in L^1([0,T], \mathrm{d}t)$.
\begin{longlist}
\item[(a)] If $\psi$ is positive and decreasing on $(0,T]$ and
$\varphi(t) = \int^t_0 \psi(s) \,\mathrm{d}s$ is $\beta$-H{\"o}lder continuous
with $\beta\in(0,1]$, then the sequence $(f_j)$ from \textup{(3.7)} satisfies
$\| f_j \|
\leq C j^{-\beta}$. If $\beta> 1/2$, then
(A1) is satisfied with $\vartheta= \beta$ and $\gamma= 0$.

\item[(b)] If $\psi(0) = 0$, $\psi$ is $\beta$-H{\"o}lder continuous with
$\beta\in(0,1]$ and $\psi$ is differentiable on $(0,T]$ such that $\psi'$
is positive and decreasing on $(0,T]$, then
(A1) is satisfied with $\vartheta= 1 + \beta$ and $\gamma= 0$.
\end{longlist}
\end{Lem}
\begin{pf}
Let $\lambda_j = (\curpi(j - 1/2)/T)^{-2}$. $(a)$ For $t
\leq
\sqrt{\lambda_j}$, we have
\[
| f_j (t) | \leq\sqrt{2/T} \varphi\bigl( \sqrt{\lambda_j}\bigr) .
\]
Using the second integral mean value formula, we obtain, for $t \in
[\sqrt
{\lambda_j}, T]$ and some\break $\delta_j \in[\sqrt{\lambda_j}, t ]$,
\begin{eqnarray*}
| f_j (t)| & \leq& \sqrt{2/T} \biggl( \biggl| \int^{\sqrt{\lambda_j}}_{0} \psi(s)
\cos\bigl(( t - s )/\sqrt{\lambda_j}\bigr) \,\mathrm{d}s \biggr| + \biggl|
\int^t_{\sqrt{\lambda_j}} \psi(s) \cos\bigl((t-s)/ \sqrt{\lambda_j} \bigr)\,\mathrm{d}s\biggl|
\biggr) \\
& = & \sqrt{2/T} \biggl( \biggl| \int^{\sqrt{\lambda_j}}_{0} \psi(s) \cos\bigl(( t
- s
)/\sqrt{\lambda_j}\bigr) \,\mathrm{d}s \biggr| + \psi\bigl(\sqrt{\lambda_j} \bigr) \biggr|
\int^{\delta_j}_{\sqrt{\lambda_j}} \cos\bigl(( t- s )/\sqrt{\lambda_j}\bigr)\,\mathrm{d}s\biggr|
\biggr) \\
& \leq& \sqrt{2/T} \bigl( \varphi\bigl( \sqrt{\lambda_j}\bigr) + 2 \sqrt{\lambda_j}
\psi\bigl( \sqrt{\lambda_j} \bigr)\bigr) \\
& \leq& 3 \sqrt{2/T} \varphi\bigl( \sqrt{\lambda_j} \bigr) .
\end{eqnarray*}
Consequently,
\[
\| f_j \| \leq3 \sqrt{2/T} \varphi\bigl(\sqrt{\lambda_j}\bigr) \leq C
\lambda
_j^{\beta/2} .
\]

$(b)$ The function $\psi$ is absolutely continuous on $[0,T]$, so
an integration by parts yields
\[
f_j (t) = \sqrt{2 \lambda_j/T} \int^t_0 \psi' (s) \sin\bigl((
t-s)/\sqrt
{\lambda_j}\bigr)\, \mathrm{d}s.
\]
Arguing as in $(a)$ (with $\psi$ replaced by $\psi'$), we deduce that
\[
\| f_j \| \leq3 \sqrt{2 \lambda_j/T} \psi\bigl( \sqrt{\lambda_j}\bigr) \leq C
\lambda^{(1+ \beta)/2}_j .
\]
\upqed
\end{pf}

Now, let $\psi(t) = t^{\rho- 1/2}$ with $\rho\in(0, \infty)$. Then
%
\begin{equation}
X_t = X_t^\rho= \int^t_0 (t-s)^{\rho- 1/2} \,\mathrm{d} W_s,\qquad t \in[0,T]
\end{equation}
so that $X^\rho$ is a \textit{Riemann--Liouville} process of order
$\rho$.
Using the
$(\rho\wedge\frac{1}{2})$-H{\"o}lder continuity of the application $t
\mapsto X^\rho_t$ from $[0,T]$ into $L^{2}(\mathbb{P})$ and the Kolmogorov
criterion, we can check that
$X^\rho$ has a pathwise continuous modification.
\begin{Lem}\label{lem4}
Let $\psi(t) = t^{\rho- 1/2}$, $\rho\in(0, \infty)$. Then the sequence
$(f_j)$ from (3.7)
satisfies (\textup{A2}) with $a = \min\{ 1, \rho+ 1/2\}$, $b = 1$ and (\textup{A1})
for $\rho\in(0, 3/2]$
with $\vartheta= \rho+ 1/2$ and $\gamma=0$.
\end{Lem}
\begin{pf}
This is an immediate consequence of Lemmas \ref{lem2} and \ref{lem3}.
\end{pf}

We deduce, for Riemann--Liouville processes of order $\rho\in(0, 3/2]$,
that
%
\begin{equation}\label{eq3.9}
e^{(\mathrm{prod})} _{ N,r} (RL)= O \biggl( \frac{ ( \log\log N)^{\rho+ 1/2} }{
(\log N)^\rho} \biggr),
\end{equation}
while for every $\rho\in(0, \infty)$ (see \cite{LI},
Graf \textit{et al.} \cite{GRAF2}),
%
\begin{equation}\label{eq3.10}
e_{ N,r} (RL) \approx( \log N)^{-\rho} .
\end{equation}

To go beyond $\rho= 3/2$, we must slightly change the way we quantize.
Let $\psi(t) = t^{\rho- 1/2}$, with $\rho> 3/2$, and choose $k \in
\mathbb{N}$
such that $k + 1/2 < \rho\leq k + 3/2$.
Set $\lambda_j = ( \curpi(j - 1/2)/T)^{-2}$. For
$k \in\{ 2 n - 1, 2 n \}$ $n \in\mathbb{N}$, integration by parts
yields the expansion
\begin{eqnarray*}
f_j (t) & = & \sum^n_{m=1} (-1)^{m-1} \lambda^m_j \sqrt{2/T} \psi
^{(2m-1)} (t) + (-1)^n \lambda^n_j \sqrt{2/T} \int^t_0 \psi^{(2n)} (s)
\cos\bigl(( t - s)/
\sqrt{\lambda_j}\bigr) \,\mathrm{d}s \\
& =: & g_j (t) + h_j(t) ,\qquad t \in[0,T] .
\end{eqnarray*}
Since $\psi^{(2n)} (t) = C t^{\beta-1}$ if $k = 2n-1$ and $\psi^{(2n)}
(t) = C t^\beta$ if $k = 2 n$ with $\beta= \rho- k - 1/2 \in(0,1]$, we
deduce from Lemma \ref{lem2} and Lemma \ref{lem3} that the sequence $(h_j)$ in $C([0,T])$
satisfies (A1) with
$\vartheta= \rho+ 1/2$, $\gamma= 0$ and (A2) with $a = \rho- k -
1/2$, $b = -k$ if $k = 2n-1$ and
$a = 1$, $b = - k + 1$ if $k = 2 n$.
Clearly, the sequence $(g_j)$ also satisfies the conditions (A1) and
(A2) (with $\vartheta= 2$, $\gamma= 0$, $b = -2$ and
$a = \rho- k - 1/2$ if $k = 2 n - 1$ and $a = 1$ if $k=2n$).
Consequently, there exist centered
continuous Gaussian processes $U = (U_t)_{t \in[0,T]}$ and $Z$ such that
$U = \sum^\infty_{j=1} \xi_j g_j $ a.s., $Z = \sum^\infty
_{j=1} \xi_j h_j $ a.s.,
%
\begin{equation}
X = X^\rho\stackrel{d}{=} U + Z
\end{equation}
and $U \in\mbox{span} \{ \psi^{(2m-1)} \dvtx m = 1 , \ldots, n \}$
a.s.
Observe that
\[
U = \sum^n_{m=1} (-1)^{m-1} \sqrt{2/T} \psi^{(2m-1)} \eta_m,
\]
where $\eta_m = \sum^\infty_{j=1} \lambda^m_j \xi_j$ is
$\mathcal{N}(0, \sum^\infty_{j=1} \lambda^{2m}_j)$-distributed.

Now use, for example, $[N^{1/2n}]$-quantizations of $\eta_m$ and a $[\sqrt
{N}]$-product quantization of $Z$ for the quantization of $X$
(which is clearly not optimal in practise, but remains rate optimal).
Let $\widehat\eta_m$ be an $L^r$-optimal $[N^{1/2n}]$-quantization
for $\eta_m$,
\[
\hat{U}^{\sqrt{N}} := \sum^n_{m=1} (-1)^{m-1} \sqrt{2/T} \psi
^{(2m-1)} \hat{\eta}_m,
\]
and let $\hat{Z}^{\sqrt{N}}$ be the $L^r$-product $[\sqrt
{N}]$-quantization of $Z$ from Theorem \ref{thm2}. A (modified)\break
$L^r$-product $N$-quantization
of $X$ with respect to $(f_j)$ is then defined by
%
\begin{equation}
\widehat X:= \hat{U}^{\sqrt{N} } + \hat{Z}^{\sqrt{N}} .
\end{equation}
Using Theorem \ref{thm2}, we can show for the quantization error,
that
\begin{eqnarray*}
\| U + Z - \widehat X\|_{L^r_E}
& \leq& C \bigl( \bigl\| U - \hat{U}^{\sqrt{N}} \bigr\|_{L^r_E} + \bigl\| Z - \hat
{Z}^{\sqrt{N}}
\bigr\|_{L^r_E} \bigr) \\
& \leq& C \Biggl( \sum^n_{m=1} \sqrt{2/T} \bigl\| \psi^{(2m-1)} \bigr\| \| \eta
_m - \hat{\eta}_m \|_{L^r} + \bigl\| Z - \hat{Z}^{\sqrt{N}} \bigr\|_{L^r_E}\Biggr) \\
& \leq&  \frac{C}{N^{1/2n}}  +
 \frac{C(\log\log\sqrt{N})^{\rho+1/2}}{(\log\sqrt
{N})^\rho}  \\
& \leq&  \frac{C ( \log\log N)^{\rho+ 1/2}}{(\log
N)^\rho}
\end{eqnarray*}
so that, with the above modification, (3.9) remains true for $\rho> 3/2$.

Now, consider the \textit{stationary Ornstein--Uhlenbeck process} as the
solution of the Langevin equation
\[
\mathrm{d} X_t = - \beta X_t \,\mathrm{d}t + \sigma \,\mathrm{d} W_t ,\qquad t
\in[0,T],
\]
with $X_0$ independent of $W$ and $\mathcal{N}(0, \frac{\sigma
^{2}}{2 \beta
}$)-distributed, $\sigma> 0$, $\beta> 0$. It admits the explicit representation
%
\begin{equation}
X_t = e^{- \beta t} X_0 + \sigma \mathrm{e}^{-\beta t} \int^t_0 \mathrm{e}^{\beta s} \,\mathrm{d} W_s
\end{equation}
and
\[
\mathbb{E}X_s X_t = \frac{\sigma^{2}}{2 \beta} \mathrm{e}^{-\beta| s-t | }.
\]
By Lemma \ref{lem2}, the admissible sequence
\[
f_0 (t) = \frac{\sigma}{\sqrt{2a}} \mathrm{e}^{-\beta t} ,\qquad f_j (t) =
\sigma\sqrt{\frac{2}{T}} \int^t_0 \mathrm{e}^{-\beta(t-s)} \cos\biggl( \frac
{\curpi (j-1/2)s}{T} \biggr)\, \mathrm{d}s, \qquad j\geq1,
\]
%
%
%
satisfies conditions (A1) and (A2) with $\vartheta= 1$, $\gamma= 0$,
$a = 1$ and $b = 1$. Consequently,
%
\begin{equation}\label{eq3.14}
e^{(\mathrm{prod})}_{ N,r} (OU) = O \biggl( \frac{ \log\log N}{(\log N)^{1/2}}
\biggr),
\end{equation}
while (see Graf \textit{et al.} \cite{GRAF2})
%
\begin{equation}\label{eq3.15}
e_{ N,r} (OU) \approx( \log N)^{-1/2} .
\end{equation}

\printhistory

\end{document}